\documentclass[12pt,a4paper,reqno]{amsart}
\usepackage{amsmath}
\usepackage{amssymb}
\usepackage{amsthm}
\usepackage{enumitem}
\usepackage{etoolbox}
\usepackage{hyperref}
\hypersetup{
  colorlinks=true,
  citecolor=blue,
  linkcolor=blue,
  urlcolor=blue
}
\numberwithin{equation}{section}
\theoremstyle{plain}
\newtheorem{theorem}{Theorem}[section]
\newtheorem{lemma}{Lemma}[section]
\newtheorem{corollary}{Corollary}[section]
\newtheorem{proposition}{Proposition}[section]
\newtheorem{example}{Example}[section]
\theoremstyle{definition}
\newtheorem*{acknowledgment}{Acknowledgments}
\newtheorem{remark}{Remark}[section]
\makeatletter
\def\@settitle{%
  \begin{center}%
    \baselineskip14\p@
    {\Large \@title\par}%
  \end{center}%
}
\patchcmd{\@setauthors}{\MakeUppercase}{}{}{}
\patchcmd{\section}{\scshape}{}{}{}
\makeatother
\begin{document}
\title
[{Restricted linear congruences over finite-norm Dedekind domains}]
{Restricted linear congruences over finite-norm Dedekind domains}

\author
[Enci Wang]
{Enci Wang}

\address{Department of Mathematics, South China University of Technology, Guangzhou, Guangdong 510640, China}
\email{maecw@mail.scut.edu.cn}

\keywords{restricted linear congruences; Ramanujan sums;
Dedekind domain; finite norm property; discrete Fourier transform}
\subjclass[2020]{11D79, 11A25, 11L03, 11R04}

\begin{abstract}
Let \(R\) be a Dedekind domain with the finite norm property, and let
\(\mathfrak{m}_1,\dots,\mathfrak{m}_k\) be pairwise coprime non-zero ideals
of \(R\). Set \(\mathfrak{m}=\mathfrak{m}_1\cdots\mathfrak{m}_k\). Given
\(a_{ij},b_i\in R\) and ideals \(\mathfrak{T}_{ij}\mid\mathfrak{m}_i\)
\((1\le i\le k,\ 1\le j\le n)\), we count the solutions
\((x_1,\dots,x_n)\in (R/\mathfrak{m})^n\) of the system
\[
\sum_{j=1}^{n}a_{ij}x_j\equiv b_i\pmod{\mathfrak{m}_i}\quad(1\le i\le k)
\]
subject to the restrictions \((x_j)+\mathfrak{m}_i=\mathfrak{T}_{ij}\).
Using generalized Ramanujan sums over ideals and the discrete Fourier
transform of arithmetic functions, we obtain an explicit formula for the
number of solutions. As a special case, we also derive an explicit product
formula for a single restricted linear congruence with arbitrary
coefficients.
\end{abstract}

\maketitle
\section{Introduction}\label{sec1}
The study of linear congruences has a long history in classical number theory. 
Lehmer \cite{Lehmer1913} established the classical solvability criterion for the linear congruence
\[
a_1x_1 + \cdots + a_kx_k \equiv b \pmod{n},\tag{1.1}\label{1.1}
\]
where \(a_1,\dots,a_k,b,n\in\mathbb{Z}\) and \(n\ge 1\). 
A solution \((x_1,\dots,x_k)\in(\mathbb{Z}/n\mathbb{Z})^k\) exists if and only if \(\ell\mid b\),
where \(\ell=(a_1,\dots,a_k,n)\); in that case the number of solutions equals \(\ell n^{k-1}\).

Restricted linear congruences, in which the variables are subject to additional arithmetic constraints, arose as a natural refinement of this problem. 
It was shown by Rademacher \cite{Rademacher1925} and Brauer \cite{Brauer1926} that the number \(N_n(k,b)\) of solutions of the congruence
\[
x_1 + \cdots + x_k \equiv b \pmod{n}\tag{1.2}\label{1.2}
\]
with the restrictions \((x_i,n)=1\) for every \(1 \leq i \leq k\) is given by
\[
N_n(k,b)=\frac{\varphi(n)^k}{n}
\prod_{\substack{p\mid n\\ p\mid b}}
\left(1-\frac{(-1)^{k-1}}{(p-1)^{k-1}}\right)
\prod_{\substack{p\mid n\\ p\nmid b}}
\left(1-\frac{(-1)^k}{(p-1)^k}\right),\tag{1.3}\label{1.3}
\]
where \(\varphi(n)\) is Euler's totient function and the products are taken over all prime divisors \(p\) of \(n\). 
This formula was later rediscovered by Dixon \cite{Dixon1960} and Rearick \cite{Rearick1963}. 
This result was expressed in terms of Ramanujan sums by Nicol and Vandiver \cite{Nicol1954} and reproved by Cohen \cite{Cohen1955}:
\[
N_n(k,b) = \frac{1}{n}\sum_{d\mid n} c_d(b)
\left(c_n\left(\frac{n}{d}\right)\right)^k,
\tag{1.4}\label{1.4}
\]
where \(c_n(m)\) denotes the classical Ramanujan sum \cite{Ramanujan1918}.  
This reformulation made the connection between restricted linear congruences and trigonometric sums explicit. 
The special case \(k=2\) was studied independently by Alder \cite{Alder1958}, Deaconescu \cite{Deaconescu2000}, and Sander \cite{Sander2009}.

The case of arbitrary coefficients,
\[
a_1x_1 + \cdots + a_kx_k \equiv b \pmod{n}, \quad (x_i, n) = t_i \quad (1 \leq i \leq k),
\tag{1.5}\label{1.5}
\]
turned out to be considerably more delicate. 
Sburlati \cite{Sburlati2003} treated this general form under certain assumptions on the prime factors of \(n\) relative to the coefficients \(a_i\) and the prescribed values \(t_i\), but the proof given there is incomplete.
Special cases were settled by Cohen \cite{Cohen1959}, who considered prime coefficients \(a_i\) with \(a_i\mid n\) and \(t_i=1\); 
by Jacobson and Williams \cite{Jacobson1972}, who treated the case \((a_1,\dots,a_k)=1\) and \(t_i=1\); 
by Sander and Sander \cite{Sander2013}, who investigated the bivariate case \(k=2\) with \(t_1=t_2=1\); 
and by Sun and Yang \cite{Sun2014}, who studied
the case \(a_i=1\) for all \(i\). 
Bibak, Kapron, Srinivasan, Tauraso, and T\'oth \cite{Bibak2017} then employed properties of Ramanujan sums together with the discrete Fourier transform (DFT) of arithmetic functions to obtain a complete and explicit formula for the number of solutions of (\ref{1.5}) for arbitrary integers \(a_i,t_i,b,n\). 
Their approach unified and generalized all of the preceding results.
It also provided necessary and sufficient conditions for the existence of solutions.

Namboothiri \cite{Namboothiri2018} pursued a generalization in a different direction, considering congruences of the form
\[
a_1^s x_1 + \cdots + a_k^s x_k \equiv b \pmod{n^s}\tag{1.6}\label{1.6}
\]
with the generalized restriction \((x_i, n^s)_s = t_i^s\), where \((a,b)_s\) denotes
the largest integer \(\ell\) such that \(\ell^s\) divides both \(a\) and \(b\). 
Using Cohen's generalized Ramanujan sums \cite{Cohen1949} and the theory of \((r,s)\)-even functions, he extended the results of \cite{Bibak2017} to this higher-power setting.

Most recently, Babu, Bera, and Sury \cite{Babu2024} studied systems of restricted linear congruences with pairwise coprime moduli,
\[
\sum_{j=1}^n a_{ij}x_j \equiv b_i \pmod{m_i},\quad
(x_j, m_i) = t_{ij}\quad (1 \leq i \leq k,\ 1 \leq j \leq n).\tag{1.7}\label{1.7}
\]
By combining the DFT approach with the Chinese remainder theorem and multivariable M\"obius inversion, they derived a general formula for the number of solutions. 
They also developed a parallel theory over the polynomial ring \(\mathbb{F}_q[t]\) over a finite field.

In the present paper we work over a Dedekind domain \(R\) with the finite norm property, that is, a Dedekind domain in which every non-zero ideal \(\mathfrak{m}\) has finite norm \(N(\mathfrak{m}) = |R/\mathfrak{m}|\).
Let \(a_{ij}, b_i \in R\), and let \(\mathfrak{m}_1, \dots, \mathfrak{m}_k\) be pairwise coprime non-zero ideals of \(R\). 
We study the system
\begin{equation}
\begin{aligned}
a_{11}x_1 + \cdots + a_{1n}x_n &\equiv b_1 \pmod{\mathfrak{m}_1}\\
&\vdots\\
a_{k1}x_1 + \cdots + a_{kn}x_n &\equiv b_k \pmod{\mathfrak{m}_k}
\end{aligned}\tag{1.8}\label{1.8}
\end{equation}
subject to the restrictions \((x_j) + \mathfrak{m}_i = \mathfrak{T}_{ij}\),
where \(\mathfrak{T}_{ij} \mid \mathfrak{m}_i\) for all \(1 \leq i \leq k\), \(1 \leq j \leq n\), and count its solutions \((x_1, \dots, x_n) \in (R/\mathfrak{m})^n\), with \(\mathfrak{m} = \mathfrak{m}_1\mathfrak{m}_2\cdots\mathfrak{m}_k\).

\begin{theorem}\label{mainthm}
Let \(\mathfrak{m}_1, \dots, \mathfrak{m}_k\) be pairwise coprime non-zero ideals of
\(R\), and set \(\mathfrak{m} = \mathfrak{m}_1 \cdots \mathfrak{m}_k\). 
Let \(\mathfrak{T}_{ij} \mid \mathfrak{m}_i\) for \(1 \leq i \leq k\), \(1 \leq j \leq n\). 
Then the number of solutions of the system
(\ref{1.8}) subject to the restrictions
\[
(x_j) + \mathfrak{m}_i = \mathfrak{T}_{ij} \quad (1 \leq i \leq k,\ 1 \leq j \leq n)
\]
is given by
\[
\frac{1}{N(\mathfrak{m})} \prod_{j=1}^{n}
\frac{\varphi\!\left(\frac{\mathfrak{m}}{\mathfrak{T}_j}\right)}
{\varphi\!\left(\frac{\mathfrak{m}}{\mathfrak{T}_j\mathfrak{D}_j}\right)}
\sum_{\mathfrak{D} \mid \mathfrak{m}} \eta(b, \mathfrak{D})
\prod_{l=1}^{n}
\eta\!\left(\frac{\mathfrak{m}}{\mathfrak{D}}, \frac{\mathfrak{m}}{\mathfrak{T}_l\mathfrak{D}_l}\right),
\]
where
\[
\mathfrak{T}_j = \prod_{i=1}^{k} \mathfrak{T}_{ij},\quad
\mathfrak{D}_j = \prod_{i=1}^{k} \mathfrak{D}_{ij},\quad
\mathfrak{D}_{ij} = (a_{ij}) + \frac{\mathfrak{m}_i}{\mathfrak{T}_{ij}}
\quad (1 \leq i \leq k,\ 1 \leq j \leq n),
\]
and \(b\) is the unique residue class modulo \(\mathfrak{m}\) determined by the congruences
\[
b \equiv b_i \pmod{\mathfrak{m}_i} \quad (1 \leq i \leq k).
\]
\end{theorem}
\begin{remark}\label{rem:spec}
Here \(\varphi(\cdot)\) denotes the Euler function of ideals and \(\eta(\cdot,\cdot)\) denotes the generalized Ramanujan sum over ideals; both are defined in detail in Section 2.
Theorem \ref{mainthm} subsumes the earlier results discussed above. 
For \(R=\mathbb{Z}\) and \(k=1\), it reduces to the main theorem of \cite{Bibak2017}; for \(R=\mathbb{Z}\) and for \(R=\mathbb{F}_q[t]\), it recovers \cite[Theorems 3 and 7]{Babu2024}, respectively. 
We emphasize that no solvability
hypothesis is imposed in the statement: the formula always counts the (possibly
empty) solution set of (\ref{1.8}).
\end{remark}
The proof of Theorem \ref{mainthm} proceeds by two reductions. 
In Section \ref{sec3} we first solve the single-variable congruence \(ax \equiv b \pmod{\mathfrak{m}}\) with the prescribed condition \((x) + \mathfrak{m} = \mathfrak{T}\), where \(\mathfrak{T}\mid\mathfrak{m}\), by a local computation via the Chinese remainder theorem (Theorem \ref{the3.1}); this accounts for the product of quotients of Euler functions appearing in the formula. 
We then use the convolution theorem for the DFT together with the inversion formula for \(\mathfrak{m}\)-even functions (Theorems \ref{the2.3} and \ref{the2.4}) to count the solutions of the diagonal congruence \(x_1 + \cdots + x_n \equiv b \pmod{\mathfrak{m}}\) with prescribed gcd ideals (Theorem \ref{the3.2}). 
Combining these two results yields Theorem \ref{mainthm}. 
As an application, Section 3 also derives an explicit product formula for the number of solutions of a single restricted congruence with arbitrary coefficients.

The paper is organized as follows. 
Section \ref{sec2} collects the necessary preliminaries on ideals of Dedekind domains: the M\"obius and Euler functions, the generalized Ramanujan sums and their orthogonality relations, and the discrete Fourier transform of \(\mathfrak{m}\)-periodic and \(\mathfrak{m}\)-even functions. 
In Section \ref{sec3} we prove Theorems \ref{the3.1} and \ref{the3.2} and complete the proof of Theorem \ref{mainthm}. 
We then derive the single-congruence application and conclude with two explicit examples, over \(\mathbb{Z}[\sqrt{-5}]\) and over the \(3\)-adic integers \(\mathbb{Z}_3\), in which the main formula is verified by direct enumeration.
\section{Notation and preliminaries}\label{sec2}
Throughout this paper, \(R\) denotes a Dedekind domain with the finite norm property. 
For a non-zero ideal \(\mathfrak{m}\) of \(R\), let
\(N(\mathfrak{m})=|R/\mathfrak{m}|\) denote its finite norm, and for a non-zero ideal \(\mathfrak{a}\) of \(R\), denote by \(v_\mathfrak{p}(\mathfrak{a})\) the exponent of the prime ideal \(\mathfrak{p}\) in the factorization of \(\mathfrak{a}\). 
For \(a\in R\setminus\{0\}\), set \(v_\mathfrak{p}(a)=v_\mathfrak{p}((a))\), with the convention \(v_\mathfrak{p}(0)=\infty\).
\subsection{Ramanujan sums}\label{sec2.1}
For any non-zero ideal \(\mathfrak{m}\) of \(R\), the M\"obius function \(\mu(\mathfrak{m})\) is defined by
\[
\mu(\mathfrak{m}) =
\begin{cases}
1, & \text{if } \mathfrak{m} = R, \\
(-1)^k, & \text{if } \mathfrak{m} = \mathfrak{p}_1\cdots\mathfrak{p}_k
\text{ for distinct prime ideals } \mathfrak{p}_1,\dots,\mathfrak{p}_k \text{ of } R, \\
0, & \text{otherwise}.
\end{cases}
\]
For any non-zero ideal \(\mathfrak{m}\) of \(R\), we have
\begin{equation}
\sum_{\mathfrak{D} \mid \mathfrak{m}} \mu(\mathfrak{D}) =
\begin{cases}
1, & \text{if } \mathfrak{m} = R, \\
0, & \text{otherwise},
\end{cases}\tag{2.1}
\end{equation}
where the sum runs over all ideal divisors \(\mathfrak{D}\) of \(\mathfrak{m}\). 
For an arithmetic function \(f\) defined on the set of ideals, set \(F(\mathfrak{m})=\sum_{\mathfrak{D}\mid \mathfrak{m}}f(\mathfrak{D})\). 
Then the M\"obius inversion formula reads
\[
f(\mathfrak{m})=\sum_{\mathfrak{D}\mid \mathfrak{m}}\mu(\mathfrak{D})\,F\!\left(\frac{\mathfrak{m}}{\mathfrak{D}}\right).\tag{2.2}\label{2.2}
\]
The Euler function \(\varphi(\mathfrak{m})\) of a non-zero ideal \(\mathfrak{m}\) is defined as the cardinality of the unit group of the quotient ring \(R/\mathfrak{m}\). 
By \cite[Theorem 1.19]{Narkiewicz2004}, we have
\[
\varphi(\mathfrak{m}) = N(\mathfrak{m}) \prod_{\mathfrak{p} \mid \mathfrak{m}} \left(1 - N(\mathfrak{p})^{-1}\right)
= \sum_{\mathfrak{D} \mid \mathfrak{m}} \mu(\mathfrak{D}) N\!\left(\frac{\mathfrak{m}}{\mathfrak{D}}\right),\tag{2.3}\label{2.3}
\]
where the product runs over all distinct prime ideals \(\mathfrak{p}\) dividing \(\mathfrak{m}\). 
By the above M\"obius inversion formula (\ref{2.2}), we obtain
\[
\sum_{\mathfrak{D}\mid\mathfrak{m}}\varphi(\mathfrak{D})=N(\mathfrak{m}).\tag{2.4}
\]

Let \(\mathfrak{m}\) be a non-zero proper ideal of \(R\). 
For any complex-valued additive character \(\psi_\mathfrak{m}\) of the quotient ring \(R/\mathfrak{m}\) and any \(a \in R\), define the character \(\psi_\mathfrak{m}^{(a)}\) by
\[
\psi_\mathfrak{m}^{(a)}(\bar{\alpha}) = \psi_\mathfrak{m}(\overline{a\alpha}) \quad \text{for all } \bar{\alpha} \in R/\mathfrak{m}.
\]
The principal character of \(R/\mathfrak{m}\) is denoted by \(\psi_0\), i.e., \(\psi_0(\bar{\alpha})=1\) for all \(\bar{\alpha} \in R/\mathfrak{m}\).
An additive character \(\psi_\mathfrak{m}\) of \(R/\mathfrak{m}\) is called primitive if, for every \(a \notin \mathfrak{m}\), the character \(\psi_\mathfrak{m}^{(a)}\) is non-principal on the additive group of \(R / \mathfrak{m}\mathfrak{n}^{-1}\), where \(\mathfrak{n} = (aR, \mathfrak{m})\). 
Since \(\mathfrak{n}\mid\mathfrak{m}\), we have \(\mathfrak{m}\subseteq\mathfrak{m}\mathfrak{n}^{-1}\), so \(R/\mathfrak{m}\mathfrak{n}^{-1}\) is a quotient of \(R/\mathfrak{m}\); moreover \(\psi_\mathfrak{m}^{(a)}\) is trivial on the kernel \(\mathfrak{m}\mathfrak{n}^{-1}/\mathfrak{m}\) because \(\beta a\in\mathfrak{m}\) for every \(\beta\in\mathfrak{m}\mathfrak{n}^{-1}\), and hence it descends to that quotient.

As shown in \cite[Proposition 1.1]{Zheng2023}, a primitive additive character of \(R/\mathfrak{m}\) exists. 
Any additive character \(\psi_\mathfrak{m}\) of \(R/\mathfrak{m}\) may be extended to a function on \(R\) by setting \(\psi_\mathfrak{m}(a) = \psi_\mathfrak{m}(\bar{a})\) for all \(a \in R\). 
For each non-zero ideal \(\mathfrak{m}\), we fix a primitive additive character of \(R/\mathfrak{m}\), extend it to a function on \(R\) as above, and denote this extension by \(e_\mathfrak{m}\). 
By construction, \(e_\mathfrak{m}\) satisfies
\[
\begin{cases}
e_\mathfrak{m}(\alpha) = e_\mathfrak{m}(\beta), & \text{if } \alpha \equiv \beta \pmod{\mathfrak{m}}, \\
e_\mathfrak{m}(\alpha + \beta) = e_\mathfrak{m}(\alpha) e_\mathfrak{m}(\beta), & \text{for all } \alpha, \beta \in R, \\
|e_\mathfrak{m}(\alpha)| = 1, & \text{for all } \alpha \in R.
\end{cases}
\]
Thus \(e_\mathfrak{m}\) is the natural analogue of the classical exponential \(e(x)=\exp(2\pi ix)\). 
The primitive characters enjoy the following fundamental property (see \cite[Corollary 2.5]{Zheng2023}).
\begin{proposition}\label{prop2.5}
Let \(\mathfrak{m}\) be a non-zero ideal of \(R\), and let \(e_{\mathfrak{m}}\) be a primitive additive character modulo \(\mathfrak{m}\). 
For every \(a \in R\), we have
\[
\sum_{\alpha\bmod \mathfrak{m}} e_{\mathfrak{m}}(a\alpha) =
\begin{cases}
N(\mathfrak{m}), & \text{if } a \equiv 0 \pmod{\mathfrak{m}},\\
0, & \text{otherwise},
\end{cases}\tag{2.5}
\]
where \(\alpha \bmod \mathfrak{m}\) means that \(\alpha\) runs through a complete system of residue classes modulo \(\mathfrak{m}\).
\end{proposition}
We shall use the following standard facts on primitive additive characters (see \cite[Lemmas 2.3 and 2.7]{zheng2022}). 
First, every additive character of \(R/\mathfrak{m}\) is of the form \(e_\mathfrak{m}^{(a)}\) for a unique \(a\bmod\mathfrak{m}\), and \(e_\mathfrak{m}^{(a)}\) is again primitive if and only if \(\bar a\) is a unit of \(R/\mathfrak{m}\). 
Second, if \(\mathfrak{D}\mid\mathfrak{m}\) and \(d\in R\) satisfies \((dR,\mathfrak{m})=\mathfrak{D}\), then the character \(\beta\mapsto e_\mathfrak{m}(d\beta)\) is primitive modulo \(\mathfrak{m}\mathfrak{D}^{-1}\).
For any \(a \in R\) and non-zero ideal \(\mathfrak{m}\) of \(R\), the Ramanujan sum \(\eta(a, \mathfrak{m})\) is defined by
\[
\eta(a, \mathfrak{m}) = \sum_{\substack{\alpha \bmod \mathfrak{m} \\ (\alpha R, \mathfrak{m}) = 1}} e_\mathfrak{m}(a\alpha) \quad \text{for } \mathfrak{m} \neq R,
\]
with the convention \(\eta(a, R) = 1\). 
This definition is independent of the choice of the primitive character
\(e_\mathfrak{m}\). 
By the first fact above, every other primitive character is \(e_\mathfrak{m}^{(u)}\) for a unit \(\bar u\in R/\mathfrak{m}\), and the substitution \(\alpha\mapsto u\alpha\) permutes the reduced residue classes modulo \(\mathfrak{m}\).
It follows directly from the change of variables \(\alpha\mapsto-\alpha\) that
\(\eta(-a,\mathfrak{D})=\eta(a,\mathfrak{D})\) for all \(a\in R\) and non-zero
ideals \(\mathfrak{D}\).

We now define the generalized Ramanujan sum indexed by pairs of ideals of \(R\). 
Let \(\mathfrak{n}\) and \(\mathfrak{m}\) be ideals of \(R\) with \(\mathfrak{m} \neq 0\), and set \(\mathfrak{D} = (\mathfrak{n}, \mathfrak{m})\). 
Since \(R/\mathfrak{m}\) is a principal ideal ring for any non-zero ideal \(\mathfrak{m}\) (see \cite[p.278]{Zariski1965}), and \(\mathfrak{D}/\mathfrak{m}\) is an ideal of \(R/\mathfrak{m}\), there exists an element \(d \in \mathfrak{D}\), called a generator of \(\mathfrak{D}/\mathfrak{m}\), whose residue class \(\bar d\) generates \(\mathfrak{D}/\mathfrak{m}\) in \(R/\mathfrak{m}\); we write \(d_\mathfrak{n}=d\). 
The following facts hold for these representatives.
\begin{lemma}\label{lem2.1}\cite[Lemma 3.3]{Zheng2023}
Let \(\mathfrak{m}\neq 0\) and \(\mathfrak{n}\) be ideals of \(R\). 
Let \(\mathfrak{D} = (\mathfrak{m}, \mathfrak{n})\), and let \(d_\mathfrak{n}\) be a generator of \(\mathfrak{D}/\mathfrak{m}\) in \(R/\mathfrak{m}\). 
Then we have
\begin{enumerate}
\item[(i)] \((d_\mathfrak{n} R, \mathfrak{m}) = (\mathfrak{n}, \mathfrak{m}) = \mathfrak{D}\).
\item[(ii)] If \(d_\mathfrak{n}'\) is another generator of \(\mathfrak{D}/\mathfrak{m}\) in \(R/\mathfrak{m}\), then \(d_\mathfrak{n}'\equiv u\,d_\mathfrak{n}\pmod{\mathfrak{m}}\) for some element \(u\in R\) with \(\bar u\) a unit of \(R/\mathfrak{m}\).
\item[(iii)] If \(\mathfrak{p} \mid \mathfrak{m}\), then we have \(\mathfrak{p} \mid d_\mathfrak{n} R\) if and only if \(\mathfrak{p} \mid \mathfrak{n}\).
\end{enumerate}
\end{lemma}
Define the generalized Ramanujan sum \(\eta(\mathfrak{n}, \mathfrak{m})\) by
\[
\eta(\mathfrak{n}, \mathfrak{m}) = \eta(d_\mathfrak{n}, \mathfrak{m}) = \sum_{\substack{\alpha \bmod \mathfrak{m} \\ (\alpha R, \mathfrak{m}) = 1}} e_\mathfrak{m}(d_\mathfrak{n} \alpha),
\]
with the convention \(\eta(\mathfrak{n}, R) = 1\) for any ideal \(\mathfrak{n}\) of \(R\).
By Lemma \ref{lem2.1}(ii), \(\eta(\mathfrak{n}, \mathfrak{m})\) is well-defined: if \(d_\mathfrak{n}'\equiv u d_\mathfrak{n}\pmod{\mathfrak{m}}\) for some \(u\) as above, then the substitution \(\alpha\mapsto u\alpha\) permutes the reduced residue classes modulo \(\mathfrak{m}\). 
In particular, if \(\mathfrak{n} = aR\) is a principal ideal of \(R\), then \(\eta(\mathfrak{n}, \mathfrak{m}) = \eta(a, \mathfrak{m})\), consistent with the earlier definition. 
Henceforth we identify the element and ideal versions of Ramanujan sums. 
We summarize several useful properties of these sums below.

\begin{theorem}\label{the2.1}\cite[Theorem 1.6]{Zheng2023}
Let \(\mathfrak{m}\neq 0\) and \(\mathfrak{n}\) be ideals of \(R\). Then
\[
\eta (\mathfrak{n},\mathfrak{m})=\sum _{\mathfrak{D} \mid (\mathfrak{m},\mathfrak{n})} N(\mathfrak{D})\, \mu \!\left( \frac{\mathfrak{m}}{\mathfrak{D}} \right)
= \varphi(\mathfrak{m})\, \mu \!\left( \frac{\mathfrak{m}}{(\mathfrak{m},\mathfrak{n})} \right) \varphi^{-1} \!\left( \frac{\mathfrak{m}}{(\mathfrak{m},\mathfrak{n})} \right).\tag{2.6}\label{2.6}
\]
\end{theorem}
In the first identity in (\ref{2.6}), if \(\mathfrak{m} \mid \mathfrak{n}\), then using the M\"obius inversion formula (\ref{2.2}), we see that the sum over \(\mathfrak{D}\mid\mathfrak{m}\) of \(\eta(\mathfrak{n},\mathfrak{D})\) is equal to \(N(\mathfrak{m})\). If \(\mathfrak{m}\nmid\mathfrak{n}\), this sum is zero by interchanging the order of summation. 
Thus, for any non-zero ideal \(\mathfrak{m}\), \(\mathfrak{n}\), we obtain the following formula:  
\[
\sum_{\mathfrak{D}\mid\mathfrak{m}}\eta(\mathfrak{n},\mathfrak{D})=\begin{cases}
N(\mathfrak{m}), & \text{if } \mathfrak{m}\mid\mathfrak{n}, \\
0, & \text{otherwise}. 
\end{cases}\tag{2.7}\label{2.7}
\]

By Theorem \ref{the2.1}, for a fixed \(\mathfrak{n}\) the function \(\mathfrak{m}\mapsto \eta(\mathfrak{n},\mathfrak{m})\) is multiplicative. 
Consequently, to evaluate \(\eta(\mathfrak{n},\mathfrak{m})\) for a given \(\mathfrak{n}\), it suffices to consider the case where \(\mathfrak{m}\) is a power of a prime ideal. 
Applying the definition of the M\"obius function together with the first equality in (\ref{2.6}), we obtain
\begin{align*}
\eta(\mathfrak{n}, \mathfrak{p}^e) = 
\begin{cases}
N(\mathfrak{p})^e - N(\mathfrak{p})^{e-1}, & \text{if } 1 \leq e \leq v_\mathfrak{p}(\mathfrak{n}), \\
-N(\mathfrak{p})^{e-1}, & \text{if } e = v_\mathfrak{p}(\mathfrak{n}) + 1, \\
1,&\text{if }e=0,\\
0, & \text{otherwise}.
\end{cases}\tag{2.8}\label{2.8}
\end{align*}

The Ramanujan sums defined above possess several key orthogonality relations, one of which is presented here.
\begin{proposition}\label{prop2.1}\cite[Proposition 2.13]{Qi2025}
For a non-zero ideal \(\mathfrak{m}\) and \(\mathfrak{D}_1 \mid \mathfrak{m},\,\mathfrak{D}_2 \mid \mathfrak{m}\), we have
\[
\sum_{\mathfrak{D} \mid \mathfrak{m}} \eta\left( \frac{\mathfrak{m}}{\mathfrak{D}}, \mathfrak{D}_1 \right) \eta\left( \frac{\mathfrak{m}}{\mathfrak{D}_2}, \mathfrak{D} \right) =
\begin{cases}
N(\mathfrak{m}), & \text{if } \mathfrak{D}_1 = \mathfrak{D}_2, \\
0, & \text{otherwise}.
\end{cases}\tag{2.9}\label{2.9}
\]
\end{proposition}

\subsection{The discrete Fourier transform}\label{sec2.2}
A function \( f: R \to \mathbb{C} \) is called  \( \mathfrak{m} \)-periodic  if \( f(a) = f(b) \) when \(a\equiv b\pmod{\mathfrak{m}}\). 
It follows from the properties of primitive characters that \(  a \mapsto \eta(a,\mathfrak{m})\) is an \(\mathfrak{m}\)-periodic function.
Let \(e_{\mathfrak{m}}\) be a primitive additive character modulo \(\mathfrak{m}\).
The discrete Fourier transform (DFT) of an \( \mathfrak{m} \)-periodic function \( f \) is the function \( \widehat{f} = \mathcal{F}(f) \) given by
\[
\widehat{f}(a) = \sum_{\alpha\bmod \mathfrak{m}}f(\alpha)\, e_\mathfrak{m}\left(-a\alpha \right) \quad (a \in R),\tag{2.10}
\]
and \(f\) is recovered from \(\widehat f\) by the inverse discrete Fourier transform (IDFT)
\[
f(a) = \frac{1}{N(\mathfrak{m})} \sum_{\alpha\bmod \mathfrak{m}} \widehat{f}(\alpha)\, e_\mathfrak{m}\left(a\alpha \right) \quad (a \in R).\tag{2.11}\label{2.11}
\]

The Cauchy convolution of two \( \mathfrak{m} \)-periodic functions \( f_1 \) and \( f_2 \) is the \( \mathfrak{m} \)-periodic function \( f_1 \otimes f_2 \) defined by
\[
(f_1 \otimes f_2)(a) = \sum_{\alpha+\beta \equiv a \bmod{\mathfrak{m}}}f_1(\alpha) f_2(\beta) \quad (a \in R).
\]
A direct computation gives
\[
(\widehat{f_1 \otimes f_2})(b)
=\sum_{\alpha\bmod\mathfrak{m}}f_1(\alpha)\,e_\mathfrak{m}(-b\alpha)
\sum_{\beta\bmod\mathfrak{m}}f_2(\beta)\,e_\mathfrak{m}(-b\beta)
=\widehat{f}_1(b)\, \widehat{f}_2(b).
\]
More generally, if \( f_1, \ldots, f_k \) are \( \mathfrak{m}\)-periodic functions, then
\[
\mathcal{F}(f_1 \otimes \cdots \otimes f_k) = \mathcal{F}(f_1) \cdots \mathcal{F}(f_k).\tag{2.12}\label{2.12}
\]

For \( \mathfrak{D} \mid \mathfrak{m} \), let \( \rho_{\mathfrak{m},\mathfrak{D}} \) be the \( \mathfrak{m} \)-periodic function defined for every \(a\in R\) by
\[
\rho_{\mathfrak{m},\mathfrak{D}}(a) = \begin{cases} 1, & \text{if } (aR, \mathfrak{m}) = \mathfrak{D}, \\ 0, & \text{otherwise}. \end{cases}
\]
The relevance of these functions lies in the connection between their discrete Fourier transform and the Ramanujan sums defined above.
\begin{theorem}\label{the2.3}
For every \( \mathfrak{D} \mid \mathfrak{m} \) and \(a\in R\),
\[
\widehat{\rho}_{\mathfrak{m},\mathfrak{D}}(a) = \eta\left(a,\frac{\mathfrak{m}}{\mathfrak{D}}\right).\tag{2.13}
\]
In particular, the Ramanujan sum \( a \mapsto \eta(a,\mathfrak{m}) \) is the DFT of the function \( a \mapsto \rho_{\mathfrak{m},R}(a) \).   
\end{theorem}
\begin{proof}
Let \(\mathfrak{m}'=\mathfrak{m}\mathfrak{D}^{-1}\), and let \(d\) be a generator of \(\mathfrak{D}/\mathfrak{m}\) in \(R/\mathfrak{m}\), so that \((dR,\mathfrak{m})=\mathfrak{D}\). 
By \cite[Lemma 3.2]{Zheng2023}, there is an additive group isomorphism
\[
\sigma: R/\mathfrak{m}' \to \mathfrak{D}/\mathfrak{m}, \quad \overline{\beta} \mapsto \overline{d\beta},
\]
where \(\overline{\beta} = \beta + \mathfrak{m}'\) and \(\overline{d\beta} = d\beta + \mathfrak{m}\). 
The map induces a bijection between the reduced residue classes modulo \(\mathfrak{m}'\) and the classes \(\alpha\bmod\mathfrak{m}\) with \((\alpha R,\mathfrak{m})=\mathfrak{D}\); hence
\[
\widehat{\rho}_{\mathfrak{m},\mathfrak{D}}(a)
=\sum_{\substack{\beta\bmod\mathfrak{m}'\\ (\beta R,\mathfrak{m}')=R}} e_\mathfrak{m}(-ad\beta).
\]
By the facts quoted above, \(e_\mathfrak{m}(d\,\cdot)\) is a primitive character modulo \(\mathfrak{m}'\), hence equals \(e_{\mathfrak{m}'}(u\,\cdot)\) for some unit \(\bar u\) of \(R/\mathfrak{m}'\). 
Substituting \(\gamma=u\beta\), which permutes the reduced residue classes, we get
\[
\widehat{\rho}_{\mathfrak{m},\mathfrak{D}}(a)
=\sum_{\substack{\gamma\bmod\mathfrak{m}'\\ (\gamma R,\mathfrak{m}')=R}} e_{\mathfrak{m}'}(-a\gamma)
=\eta(-a,\mathfrak{m}')=\eta(a,\mathfrak{m}').
\]
This completes the proof.
\end{proof}
An \(\mathfrak{m}\)-periodic function \(f: R \to \mathbb{C}\) is called \(\mathfrak{m}\)-even if \(f(a)=f(b)\) whenever \((aR,\mathfrak{m})=(bR,\mathfrak{m})\); equivalently, \(f(a)\) depends only on the ideal \((aR,\mathfrak{m})\). 
Thus, if \(d_a\) is a generator of \((aR,\mathfrak{m})/\mathfrak{m}\) in \(R/\mathfrak{m}\), then \(f(a)=f(d_a)\). 
By Lemma \ref{lem2.1} and Theorem \ref{the2.1}, the function \( a \mapsto \eta(a,\mathfrak{m}) \) is an example of an \( \mathfrak{m}\)-even function.

For each divisor \(\mathfrak{D}\mid\mathfrak{m}\), fix a generator \(d_{\mathfrak{D}}\) of \(\mathfrak{D}/\mathfrak{m}\) in \(R/\mathfrak{m}\), so that \((d_{\mathfrak{D}}R,\mathfrak{m})=\mathfrak{D}\). 
Then for any \(a\in R\), with \(\mathfrak{D}=(aR,\mathfrak{m})\), we have \(f(a)=f(d_{\mathfrak{D}})\), hence
\[
f(a)=\sum_{\mathfrak{D}\mid\mathfrak{m}} f(d_{\mathfrak{D}})\,\rho_{\mathfrak{m},\mathfrak{D}}(a).\tag{2.14}\label{2.14}
\]
Furthermore, we have the following theorem.
\begin{theorem}\label{the2.4}
If \(f\) is an \(\mathfrak{m}\)-even function, then
\[
\widehat{f}(a) = \sum_{\mathfrak{D} \mid \mathfrak{m}} f(d_\mathfrak{D})\eta\left(a, \frac{\mathfrak{m}}{\mathfrak{D}}\right)\quad (a\in R),\tag{2.15}\label{2.15}
\]
and the IDFT is given by
\[
f(a) = \frac{1}{N(\mathfrak{m})}\sum_{\mathfrak{D} \mid \mathfrak{m}} \widehat{f}(d_\mathfrak{D})\eta\left(a, \frac{\mathfrak{m}}{\mathfrak{D}}\right)\quad (a\in R).\tag{2.16}\label{2.16}
\]
\end{theorem}
\begin{proof}
Taking the discrete Fourier transform on both sides of (\ref{2.14}) and applying Theorem \ref{the2.3}, we obtain
\[
\widehat{f}(a)=\sum_{\mathfrak{D}\mid\mathfrak{m}}f(d_{\mathfrak{D}})\,\widehat{\rho}_{\mathfrak{m},\mathfrak{D}}(a)
 =\sum_{\mathfrak{D}\mid\mathfrak{m}} f(d_{\mathfrak{D}})\,\eta\left(a,\frac{\mathfrak{m}}{\mathfrak{D}}\right),
\]
which is (\ref{2.15}).

It remains to show that \(\widehat f\) is again \(\mathfrak{m}\)-even. 
The function \(\eta(\cdot,\mathfrak{D})\) is \(\mathfrak{m}\)-periodic for any \(\mathfrak{D}\mid\mathfrak{m}\), and by Theorem \ref{the2.1} the value \(\eta(a,\mathfrak{D})\) depends only on the gcd ideal \((aR,\mathfrak{D})\). 
Let \(d_a\) be any generator of \((aR,\mathfrak{m})/\mathfrak{m}\) in \(R/\mathfrak{m}\). 
Since \(\mathfrak{D}\mid\mathfrak{m}\), we have
\[
(aR,\mathfrak{D}) = ((aR,\mathfrak{m}),\mathfrak{D}) = ((d_aR,\mathfrak{m}),\mathfrak{D}) = (d_aR,\mathfrak{D}).
\]
Accordingly, \(\eta(a,\mathfrak{D})=\eta(d_a,\mathfrak{D})\), which shows that \(\eta(\cdot,\mathfrak{D})\) is \(\mathfrak{m}\)-even. 
Consequently, \(\widehat f\) is also \(\mathfrak{m}\)-even, being a linear combination of \(\mathfrak{m}\)-even functions. 
Expanding \(\widehat f\) as in (\ref{2.14}) and applying the inverse discrete Fourier transform (\ref{2.11}), we compute
\begin{align*}
f(a)
&=\frac{1}{N(\mathfrak{m})}\sum_{\alpha\bmod\mathfrak{m}}\widehat{f}(\alpha)\,e_{\mathfrak{m}}(a\alpha)\\
&=\frac{1}{N(\mathfrak{m})}\sum_{\mathfrak{D}\mid\mathfrak{m}} \widehat{f}(d_{\mathfrak{D}})
\sum_{\alpha\bmod\mathfrak{m}}\rho_{\mathfrak{m},\mathfrak{D}}(\alpha)\,e_{\mathfrak{m}}(a\alpha)\\
&=\frac{1}{N(\mathfrak{m})}\sum_{\mathfrak{D}\mid\mathfrak{m}} \widehat{f}(d_{\mathfrak{D}})\,\eta\left(a,\frac{\mathfrak{m}}{\mathfrak{D}}\right).
\end{align*}
The last equality uses Theorem \ref{the2.3} applied with the sign reversed, together with the identity \(\eta(-a,\mathfrak{D})=\eta(a,\mathfrak{D})\). 
This proves (\ref{2.16}).
\end{proof}

\section{Restricted linear congruences over \(R\)}\label{sec3}
In this section, we employ the Ramanujan sums and the discrete Fourier transform of arithmetic functions developed in Section 2 to derive an explicit formula for the number of solutions of the system (\ref{1.8}).

We begin with the case of a single variable, which will serve as a building block for the system of linear congruences. 
The following theorem generalizes the result of \cite[Theorem 3.1]{Bibak2017} and \cite[Theorem 5]{Babu2024}.

\begin{theorem}\label{the3.1}
Let \(\mathfrak{m} \subseteq R\) be a non-zero ideal, and let \(a, b \in R\). The congruence \(ax \equiv b \pmod{\mathfrak{m}}\) has a solution \(x \in R\) satisfying \((x) + \mathfrak{m} = \mathfrak{T}\) (where \(\mathfrak{T}\) is a divisor ideal of \(\mathfrak{m}\)) if and only if \(\mathfrak{T} \supseteq (b) + \mathfrak{m}\) and \((a)+ \frac{\mathfrak{m}}{\mathfrak{T}} = \frac{(b)}{\mathfrak{T}} + \frac{\mathfrak{m}}{\mathfrak{T}}\). 
Moreover, if these conditions are satisfied, then the number of solutions modulo \(\mathfrak{m}\) is
\[
\frac{\varphi\left(\frac{\mathfrak{m}}{\mathfrak{T}}\right)}{\varphi\left(\frac{\mathfrak{m}}{\mathfrak{T}\mathfrak{D}}\right)}=N(\mathfrak{D})\prod_{\substack{\mathfrak{p}\mid \mathfrak{D}\\ \mathfrak{p}\nmid \frac{\mathfrak{m}}{\mathfrak{T}\mathfrak{D}}}}\left(1-\frac{1}{N(\mathfrak{p})}\right),
\]
where the product is over prime ideals \(\mathfrak{p}\) and \(\mathfrak{D} = (a)+ \frac{\mathfrak{m}}{\mathfrak{T}} = \frac{(b)}{\mathfrak{T}} + \frac{\mathfrak{m}}{\mathfrak{T}}\).
\end{theorem}
\begin{proof}
Suppose there exists \(x\in R\) such that \(ax\equiv b\pmod{\mathfrak{m}}\) and \((x)+\mathfrak{m}=\mathfrak{T}\). 
Then \((x)\subseteq\mathfrak{T}\) and \(\frac{(x)}{\mathfrak{T}}+\frac{\mathfrak{m}}{\mathfrak{T}}=R\). 
The congruence \(ax\equiv b\pmod{\mathfrak{m}}\) gives \((b)+\mathfrak{m}=(ax)+\mathfrak{m}\subseteq\mathfrak{T}\), which is the first condition, and
\[
\frac{(b)}{\mathfrak{T}}+\frac{\mathfrak{m}}{\mathfrak{T}}
=\frac{(ax)}{\mathfrak{T}}+\frac{\mathfrak{m}}{\mathfrak{T}}
=(a)+\frac{\mathfrak{m}}{\mathfrak{T}},
\]
where the last equality follows from \(\frac{(x)}{\mathfrak{T}}+\frac{\mathfrak{m}}{\mathfrak{T}}=R\). 
This proves the second condition.

Conversely, assume the two conditions. Set \(A=R/\mathfrak{m}\), a finite principal ideal ring, and let \(\bar{u}\in A\) generate the ideal \(\mathfrak{T}/\mathfrak{m}\). 
Multiplying the second condition by \(\mathfrak{T}\) and passing to \(A\), we obtain
\[
(\bar{a})(\bar{u})=(\bar{b}) \quad\text{in } A.
\]
Since \(A\) is a finite product of local principal ideal rings, any two generators of the same principal ideal of \(A\) differ by a unit. 
Thus there exists a unit \(\bar{v}\in A^{\times}\) such that \(\bar{a}\bar{u}\bar{v}=\bar{b}\). 
Define \(\bar{x}=\bar{u}\bar{v}\); then \(\bar{a}\bar{x}=\bar{b}\) and \((\bar{x})=(\bar{u})=\mathfrak{T}/\mathfrak{m}\).
Lifting \(\bar{x}\) to any \(x\in R\) gives a solution with the required property
\((x)+\mathfrak{m}=\mathfrak{T}\).

\medskip
Now we determine the number of solutions modulo \(\mathfrak{m}\). By the Chinese remainder theorem, the ring \(R/\mathfrak{m}\) decomposes as a product over the prime powers dividing \(\mathfrak{m}\):
\[
R/\mathfrak{m}\cong\prod_{\mathfrak{p}\mid\mathfrak{m}}
R/\mathfrak{p}^{r_{\mathfrak{p}}},\quad r_{\mathfrak{p}}=v_{\mathfrak{p}}(\mathfrak{m}).
\]
For each prime \(\mathfrak{p}\mid\mathfrak{m}\), let \(t_{\mathfrak{p}}=v_{\mathfrak{p}}(\mathfrak{T})\), \(k_{\mathfrak{p}}=r_{\mathfrak{p}}-t_{\mathfrak{p}}\), and \(d_{\mathfrak{p}}=v_{\mathfrak{p}}(\mathfrak{D})\). 
The two conditions give
\[
d_{\mathfrak{p}}=\min(v_{\mathfrak{p}}(a),k_{\mathfrak{p}})
=\min(v_{\mathfrak{p}}(b)-t_{\mathfrak{p}},k_{\mathfrak{p}}),
\]
where, as usual, \(\infty-t_{\mathfrak{p}}=\infty\).

Fix a prime \(\mathfrak{p}\mid\mathfrak{m}\), and put \(\alpha=v_{\mathfrak{p}}(a)\), \(\beta=v_{\mathfrak{p}}(b)\). 
In the discrete valuation ring \(R_{\mathfrak{p}}\) with uniformizer \(\pi\) and residue field of size \(q=N(\mathfrak{p})\), we count the residues \(x\bmod\mathfrak{p}^{r_{\mathfrak{p}}}\) such that
\[
ax\equiv b\pmod{\mathfrak{p}^{r_{\mathfrak{p}}}},\quad
\min(v_{\mathfrak{p}}(x),r_{\mathfrak{p}})=t_{\mathfrak{p}}.
\]
The first condition ensures that \(b\in\pi^{t_{\mathfrak{p}}}R_{\mathfrak{p}}\), so that \(\beta\ge t_{\mathfrak{p}}\) or \(\beta=\infty\). 
We distinguish two cases.
\begin{itemize}
\item If \(t_{\mathfrak{p}}=r_{\mathfrak{p}}\), then
\(k_{\mathfrak{p}}=d_{\mathfrak{p}}=0\), and the condition \(\min(v_{\mathfrak{p}}(x),r_{\mathfrak{p}})=r_{\mathfrak{p}}\) forces \(x\equiv0\pmod{\mathfrak{p}^{r_{\mathfrak{p}}}}\). 
This residue satisfies the congruence because \(b\in\mathfrak{p}^{r_{\mathfrak{p}}}\), and it is the only one.
\item If \(t_{\mathfrak{p}}<r_{\mathfrak{p}}\), then \(\min(v_{\mathfrak{p}}(x),r_{\mathfrak{p}})=t_{\mathfrak{p}}\) is equivalent to \(v_{\mathfrak{p}}(x)=t_{\mathfrak{p}}\). 
Write \(x=\pi^{t_{\mathfrak{p}}}y\) with \(v_{\mathfrak{p}}(y)=0\). 
Substituting into the congruence and dividing by \(\pi^{t_{\mathfrak{p}}}\) gives
\[
\pi^{\alpha}y\equiv b\pi^{-t_{\mathfrak{p}}}\pmod{\pi^{k_{\mathfrak{p}}}},
\]
where we use the convention \(\pi^{\infty}=0\).
\begin{itemize}
\item If \(d_{\mathfrak{p}}<k_{\mathfrak{p}}\), then \(\alpha=\beta-t_{\mathfrak{p}}=d_{\mathfrak{p}}\); in particular, \(\alpha,\beta<\infty\), and writing \(a=u_a\pi^{\alpha}\), \(b=u_b\pi^{\beta}\) with units \(u_a,u_b\), the congruence reduces to
\[
y\equiv u\pmod{\pi^{k_{\mathfrak{p}}-d_{\mathfrak{p}}}},\qquad
u=u_a^{-1}u_b\in R_{\mathfrak{p}}^{\times}.
\]
It has \(q^{d_{\mathfrak{p}}}\) solutions, all satisfying \(v_{\mathfrak{p}}(y)=0\) because \(u\) is a unit.
\item If \(d_{\mathfrak{p}}=k_{\mathfrak{p}}\), then \(\alpha,\beta-t_{\mathfrak{p}}\ge k_{\mathfrak{p}}\), so both sides of the congruence vanish modulo \(\pi^{k_{\mathfrak{p}}}\), and it holds for every \(y\). 
The condition \(v_{\mathfrak{p}}(y)=0\) leaves exactly \(\varphi(\mathfrak{p}^{k_{\mathfrak{p}}})=q^{k_{\mathfrak{p}}}(1-q^{-1})\) choices.
\end{itemize}
\end{itemize}
Thus, in all cases, the local number of solutions is
\(\varphi(\mathfrak{p}^{k_{\mathfrak{p}}})/\varphi(\mathfrak{p}^{k_{\mathfrak{p}}-d_{\mathfrak{p}}})\).
By the Chinese remainder theorem, the total number of solutions modulo \(\mathfrak{m}\) is therefore
\[
N=\prod_{\mathfrak{p}\mid\mathfrak{m}}
\frac{\varphi(\mathfrak{p}^{k_{\mathfrak{p}}})}
{\varphi(\mathfrak{p}^{k_{\mathfrak{p}}-d_{\mathfrak{p}}})}
=\frac{\varphi(\mathfrak{m}/\mathfrak{T})}
{\varphi(\mathfrak{m}/(\mathfrak{T}\mathfrak{D}))},
\]
where the last equality follows from (\ref{2.3}) and the multiplicativity of \(\varphi\), since \(v_{\mathfrak{p}}(\mathfrak{m}/\mathfrak{T})=k_{\mathfrak{p}}\) and \(v_{\mathfrak{p}}(\mathfrak{m}/(\mathfrak{T}\mathfrak{D}))=k_{\mathfrak{p}}-d_{\mathfrak{p}}\).
Only primes dividing \(\mathfrak{D}\) contribute nontrivially to the product,
and for such a prime the local factor equals \(N(\mathfrak{p}^{d_{\mathfrak{p}}})\) if \(d_{\mathfrak{p}}<k_{\mathfrak{p}}\), and \(N(\mathfrak{p}^{d_{\mathfrak{p}}})(1-N(\mathfrak{p})^{-1})\) if \(d_{\mathfrak{p}}=k_{\mathfrak{p}}\). 
Hence
\[
N=N(\mathfrak{D})
\prod_{\substack{\mathfrak{p}\mid\mathfrak{D}\\
\mathfrak{p}\nmid\mathfrak{m}/(\mathfrak{T}\mathfrak{D})}}
\left(1-\frac{1}{N(\mathfrak{p})}\right),
\]
which completes the proof of the formula.
\end{proof}
As a first consequence of Theorem \ref{the3.1}, we characterize when the solution, if it exists, is unique.
\begin{corollary}\label{cor3.1}
Assume that the congruence \(ax \equiv b \pmod{\mathfrak{m}}\) has a solution \(x\) with \((x)+\mathfrak{m} = \mathfrak{T}\). 
Then the number of such solutions modulo \(\mathfrak{m}\) is \(1\) if and only if one of the following two cases holds:
\begin{enumerate}
\item[(i)]\(\mathfrak{D} = R\);
\item[(ii)]\(\mathfrak{D} \) is a product of distinct prime ideals each of norm \(2\), and for every such prime \(\mathfrak{p}\mid\mathfrak{D}\) we have
\(v_{\mathfrak{p}}(\mathfrak{m}) = v_{\mathfrak{p}}(\mathfrak{T})+1\).
\end{enumerate}
\end{corollary}

\begin{proof}
Suppose a solution exists, so that the formula of Theorem \ref{the3.1} applies.
If \(\mathfrak{D}=R\), then this formula gives \(N=1\). Now we suppose \(\mathfrak{D}\neq R\). 
Then there exists at least one prime ideal \(\mathfrak{p}\) dividing \(\mathfrak{D}\), and for any such prime ideal we have \(N(\mathfrak{p})\ge 2\). 
Hence
\begin{align*}
N(\mathfrak{D})\prod_{\substack{\mathfrak{p}\mid \mathfrak{D}\\ \mathfrak{p}\nmid \frac{\mathfrak{m}}{\mathfrak{T}\mathfrak{D}}}}\left(1-\frac{1}{N(\mathfrak{p})}\right)
&\ge N(\mathfrak{D})\prod_{\mathfrak{p}\mid\mathfrak{D}}\left(1-\frac{1}{N(\mathfrak{p})}\right) \\
&\ge N(\mathfrak{D})\prod_{\mathfrak{p}\mid\mathfrak{D}}\frac{1}{N(\mathfrak{p})}
= \prod_{\mathfrak{p}\mid\mathfrak{D}} N(\mathfrak{p})^{v_{\mathfrak{p}}(\mathfrak{D})-1} \ge 1.
\end{align*}
If the number of solutions equals \(1\), then all inequalities become equalities. 
This forces \(v_{\mathfrak{p}}(\mathfrak{D})=1\) and \(N(\mathfrak{p})=2\) for each \(\mathfrak{p}\mid\mathfrak{D}\), together with \(\mathfrak{p}\nmid\frac{\mathfrak{m}}{\mathfrak{T}\mathfrak{D}}\) for each \(\mathfrak{p}\mid\mathfrak{D}\), i.e.\ \(v_{\mathfrak{p}}(\mathfrak{m})=v_{\mathfrak{p}}(\mathfrak{T})+1\).

Conversely, if either of the two conditions holds, then the formula gives the number of solutions equal to \(1\).
\end{proof}
Having settled the one-variable case, we turn to the second ingredient needed for the main theorem: the congruence \(x_1+\cdots+x_n \equiv b\pmod{\mathfrak{m}}\) with \((x_i)+\mathfrak{m}=\mathfrak{m}_i\), where \(\mathfrak{m}_i\mid\mathfrak{m}\) are given.

\begin{theorem}\label{the3.2}
Let \(\mathfrak{m} \subseteq R\) be a non-zero ideal, and let \(\mathfrak{m}_i \mid \mathfrak{m}\) for \(1 \leq i \leq n\). 
Then for any \(b \in R\), the number of solutions modulo \(\mathfrak{m}\) of the linear congruence
\[
x_1 + \cdots + x_n \equiv b \pmod{\mathfrak{m}}, \quad \text{with } (x_i) + \mathfrak{m} = \mathfrak{m}_i \ (1 \leq i \leq n),
\]
is given by
\[
\frac{1}{N(\mathfrak{m})} \sum_{\mathfrak{D} \mid \mathfrak{m}} \eta(b, \mathfrak{D}) \prod_{i=1}^{n} \eta\left(\frac{\mathfrak{m}}{\mathfrak{D}},\frac{\mathfrak{m}}{\mathfrak{m}_i}\right).
\]
\end{theorem}
\begin{proof}
The number of solutions admits the convolution representation
\[
N_\mathfrak{m}(b;\mathfrak{m}_1,\dots,\mathfrak{m}_n)
=(\rho_{\mathfrak{m},\mathfrak{m}_1}\otimes\cdots\otimes\rho_{\mathfrak{m},\mathfrak{m}_n})(b)
=\sum_{\substack{x_1,\dots,x_n \bmod{\mathfrak{m}}\\
x_1+\cdots+x_n \equiv b \pmod{\mathfrak{m}}\\
(x_i)+\mathfrak{m}=\mathfrak{m}_i,\ 1\le i\le n}}1.
\]
By the convolution theorem (\ref{2.12}) and Theorem \ref{the2.3}, its DFT is
\[
\widehat{N_\mathfrak{m}}(b;\mathfrak{m}_1,\dots,\mathfrak{m}_n)
=\prod_{i=1}^{n}\widehat{\rho}_{\mathfrak{m},\mathfrak{m}_i}(b)
=\prod_{i=1}^{n}\eta\left(b,\frac{\mathfrak{m}}{\mathfrak{m}_i}\right).
\]
Since each function \(a\mapsto\eta(a,\mathfrak{D})\) is \(\mathfrak{m}\)-even, so is \(\widehat{N_\mathfrak{m}}\). 
Applying Theorem \ref{the2.4}, we obtain
\[
N_\mathfrak{m}(b;\mathfrak{m}_1,\dots,\mathfrak{m}_n)
=\frac{1}{N(\mathfrak{m})}\sum_{\mathfrak{D}\mid\mathfrak{m}}
\widehat{N_\mathfrak{m}}(d_\mathfrak{D})\,\eta\left(b,\frac{\mathfrak{m}}{\mathfrak{D}}\right)
=\frac{1}{N(\mathfrak{m})}\sum_{\mathfrak{D}\mid\mathfrak{m}}
\eta\left(b,\frac{\mathfrak{m}}{\mathfrak{D}}\right)
\prod_{i=1}^{n}\eta\left(d_\mathfrak{D},\frac{\mathfrak{m}}{\mathfrak{m}_i}\right).
\]
By Lemma \ref{lem2.1}(i), \((d_\mathfrak{D}R,\mathfrak{m})=\mathfrak{D}\), hence \((d_\mathfrak{D}R,\mathfrak{m}/\mathfrak{m}_i)=(\mathfrak{D},\mathfrak{m}/\mathfrak{m}_i)\), so that
\[
\eta\left(d_\mathfrak{D},\frac{\mathfrak{m}}{\mathfrak{m}_i}\right)
=\eta\left(\mathfrak{D},\frac{\mathfrak{m}}{\mathfrak{m}_i}\right).
\]
Hence the formula becomes
\[
N_\mathfrak{m}(b;\mathfrak{m}_1,\dots,\mathfrak{m}_n)
=\frac{1}{N(\mathfrak{m})}\sum_{\mathfrak{D}\mid\mathfrak{m}}\eta(b,\mathfrak{D})
\prod_{i=1}^{n}\eta\left(\frac{\mathfrak{m}}{\mathfrak{D}},\frac{\mathfrak{m}}{\mathfrak{m}_i}\right),
\]
which completes the proof.
\end{proof}

We are now in a position to prove the main theorem announced in the introduction. 
The proof combines Theorem \ref{the3.1} and Theorem \ref{the3.2} through the Chinese remainder theorem.

\begin{proof}[\textbf{Proof of Theorem~\ref{mainthm}}]
Let \((x_1,\dots,x_n)\) be a solution of (\ref{1.8}) with
\((x_j)+\mathfrak{m}_i=\mathfrak{T}_{ij}\) for all \(i,j\). 
Since the \(\mathfrak{m}_i\) are pairwise coprime, for each \(j\) there is a unique residue \(y_j\bmod\mathfrak{m}\) satisfying \(y_j\equiv a_{ij}x_j\pmod{\mathfrak{m}_i}\) for all \(i\). 
Then
\[
(y_j)+\mathfrak{m}_i=(a_{ij}x_j)+\mathfrak{m}_i
=(a_{ij})\mathfrak{T}_{ij}+\mathfrak{m}_i=\mathfrak{T}_{ij}\mathfrak{D}_{ij},
\]
with \(\mathfrak{D}_{ij}=(a_{ij})+\mathfrak{m}_i/\mathfrak{T}_{ij}\) as in the
statement, and \(\sum_{j=1}^{n}y_j\equiv b_i\pmod{\mathfrak{m}_i}\) for all \(i\).

Conversely, suppose that the \(y_j\in R\) satisfy \((y_j)+\mathfrak{m}_i=\mathfrak{T}_{ij}\mathfrak{D}_{ij}\) for all \(i,j\).
Then \(\mathfrak{T}_{ij}\supseteq(y_j)+\mathfrak{m}_i\) and
\[
\frac{(y_j)}{\mathfrak{T}_{ij}}+\frac{\mathfrak{m}_i}{\mathfrak{T}_{ij}}
=\mathfrak{D}_{ij}=(a_{ij})+\frac{\mathfrak{m}_i}{\mathfrak{T}_{ij}},
\]
so, by Theorem \ref{the3.1}, the congruence \(a_{ij}x_j\equiv y_j\pmod{\mathfrak{m}_i}\) has a solution \(x_j\) with \((x_j)+\mathfrak{m}_i=\mathfrak{T}_{ij}\), and by the coprimality of the \(\mathfrak{m}_i\) one and the same \(x_j\) works for all \(i\). 
Therefore the solutions of (\ref{1.8}) are in bijection, up to the lifting count below, with the tuples \((y_1,\dots,y_n)\) satisfying
\[
\sum_{j=1}^{n}y_j\equiv b_i\pmod{\mathfrak{m}_i}\ (1\le i\le k),
\quad (y_j)+\mathfrak{m}_i=\mathfrak{T}_{ij}\mathfrak{D}_{ij}\ (1\le i\le k,\ 1\le j\le n).
\]
By the Chinese remainder theorem, the first group of congruences is equivalent to the single congruence
\[
\sum_{j=1}^{n}y_j\equiv b\pmod{\mathfrak{m}},
\]
and \((y_j)+\mathfrak{m}_i=\mathfrak{T}_{ij}\mathfrak{D}_{ij}\) for all \(i\) is equivalent to \((y_j)+\mathfrak{m}=\mathfrak{T}_j\mathfrak{D}_j\). 
By Theorem \ref{the3.2}, the number of such tuples \((y_1,\dots,y_n)\) is
\[
\frac{1}{N(\mathfrak{m})}\sum_{\mathfrak{D}\mid\mathfrak{m}}\eta(b,\mathfrak{D})
\prod_{l=1}^{n}\eta\left(\frac{\mathfrak{m}}{\mathfrak{D}},\frac{\mathfrak{m}}{\mathfrak{T}_l\mathfrak{D}_l}\right).
\]
It remains to lift each such tuple. 
For fixed \(j\), the number of \(x_j\) solving \(a_{ij}x_j\equiv y_j\pmod{\mathfrak{m}_i}\) with \((x_j)+\mathfrak{m}_i=\mathfrak{T}_{ij}\) for all \(i\) is, by Theorem \ref{the3.1} and the coprimality of the \(\mathfrak{m}_i\),
\[
\prod_{i=1}^{k}\frac{\varphi(\mathfrak{m}_i/\mathfrak{T}_{ij})}{\varphi(\mathfrak{m}_i/(\mathfrak{T}_{ij}\mathfrak{D}_{ij}))}
=\frac{\varphi(\mathfrak{m}/\mathfrak{T}_j)}{\varphi(\mathfrak{m}/(\mathfrak{T}_j\mathfrak{D}_j))},
\]
where the equality follows from the multiplicativity of \(\varphi\) for pairwise coprime ideals. 
Since the choices of the \(x_j\) for different \(j\) are independent, the number of tuples \((x_1,\dots,x_n)\) lifting a fixed tuple \((y_1,\dots,y_n)\) is
\[
\prod_{j=1}^{n}\frac{\varphi(\mathfrak{m}/\mathfrak{T}_j)}{\varphi(\mathfrak{m}/(\mathfrak{T}_j\mathfrak{D}_j))}.
\]
Multiplying the two counts, the total number of solutions of (\ref{1.8}) is
\[
\frac{1}{N(\mathfrak{m})}\prod_{j=1}^{n}\frac{\varphi(\mathfrak{m}/\mathfrak{T}_j)}{\varphi(\mathfrak{m}/(\mathfrak{T}_j\mathfrak{D}_j))}
\sum_{\mathfrak{D}\mid\mathfrak{m}}\eta(b,\mathfrak{D})
\prod_{l=1}^{n}\eta\left(\frac{\mathfrak{m}}{\mathfrak{D}},\frac{\mathfrak{m}}{\mathfrak{T}_l\mathfrak{D}_l}\right),
\]
as claimed.
\end{proof}
\begin{remark}\label{rem:orth}
Before leaving the general setting, we point out a useful consequence of the formula just proved. 
A comparison of Theorems \ref{the3.1} and \ref{mainthm} in the single-variable case yields the following orthogonality relation. 
Let \(\mathfrak{m}\) be a non-zero ideal of \(R\), \(b\in R\), and \(\mathfrak{n}\mid\mathfrak{m}\). 
The congruence \(x\equiv b\pmod{\mathfrak{m}}\) has a solution with \((x)+\mathfrak{m}=\mathfrak{n}\) if and only if \((b)+\mathfrak{m}=\mathfrak{n}\), and then the solution is unique. 
Applying Theorem \ref{mainthm} with \(a=1\) and \(\mathfrak{T}=\mathfrak{n}\), the same number equals
\[
\frac{1}{N(\mathfrak{m})}\sum_{\mathfrak{D}\mid\mathfrak{m}}
\eta(b,\mathfrak{D})\,\eta\left(\frac{\mathfrak{m}}{\mathfrak{D}},
\frac{\mathfrak{m}}{\mathfrak{n}}\right).
\]
Comparing the two expressions, we obtain
\[
\sum_{\mathfrak{D} \mid \mathfrak{m}} \eta(b, \mathfrak{D})\,
\eta\left( \frac{\mathfrak{m}}{\mathfrak{D}}, \frac{\mathfrak{m}}{\mathfrak{n}} \right)
=
\begin{cases}
N(\mathfrak{m}), & \text{if } (b) + \mathfrak{m} = \mathfrak{n}, \\
0, & \text{otherwise}.
\end{cases}
\]
For \(\mathfrak{n}=\mathfrak{m}\) this identity reduces to (\ref{2.7}).
Moreover, fixing \(\mathfrak{n}=\mathfrak{D}_2\) and summing the identity over all \(b\bmod\mathfrak{m}\) with \((b)+\mathfrak{m}=\mathfrak{D}_1\), then dividing by
\(\varphi(\mathfrak{m}/\mathfrak{D}_1)\), we obtain
\[
\sum_{\mathfrak{D} \mid \mathfrak{m}} \eta(\mathfrak{D}_1, \mathfrak{D})\,
\eta\!\left( \frac{\mathfrak{m}}{\mathfrak{D}},
\frac{\mathfrak{m}}{\mathfrak{D}_2} \right)
=
\begin{cases}
N(\mathfrak{m}), & \text{if } \mathfrak{D}_1 = \mathfrak{D}_2, \\
0, & \text{otherwise}.
\end{cases}
\]
This identity is similar to the orthogonality relation (\ref{2.9}) of Proposition \ref{prop2.1}, with the roles of the two arguments of the Ramanujan sum interchanged, but it is not identical with (\ref{2.9}), because \(\eta(\cdot,\cdot)\) is not symmetric in its two arguments. 
It also shows, in the single-variable case, how the formula vanishes when no solution exists, as anticipated in Remark \ref{rem:spec}.
\end{remark}
As an application, we specialize Theorem \ref{mainthm} to the case \(k=1\), that is, to a single restricted linear congruence
\[
a_1x_1 + \cdots + a_nx_n \equiv b \pmod{\mathfrak{m}},\tag{3.1}\label{3.1}
\]
with the restrictions \((x_i) + \mathfrak{m} = \mathfrak{T}_i\), where
\(\mathfrak{T}_i \mid \mathfrak{m}\) for all \(1 \leq i \leq n\).

If \(a_i = 0\) for all \(1 \leq i \leq n\), then a solution \((x_1, \dots, x_n)\) exists if and only if \(b \equiv 0 \pmod{\mathfrak{m}}\), and in that case the number of solutions is \(\prod_{i=1}^{n}\varphi\left(\frac{\mathfrak{m}}{\mathfrak{T}_i}\right)\).

We now assume that there exists some index \(i_0\) such that \(a_{i_0} \neq 0\).
For each prime ideal \(\mathfrak{p} \mid \mathfrak{m}\), define
\[
m_{\mathfrak{p}} = \min\left\{ v_{\mathfrak{p}}(a_i \mathfrak{T}_i) + 1 \mid 1 \leq i \leq n \right\}.
\]
Since \(a_{i_0}\neq0\), we have \(v_{\mathfrak{p}}(a_{i_0}\mathfrak{T}_{i_0})<\infty\), so \(m_{\mathfrak{p}}\) is finite for every \(\mathfrak{p}\mid\mathfrak{m}\).
Furthermore, let
\[
e_{\mathfrak{p}} = \#\left\{ i : 1 \leq i \leq n,\ v_{\mathfrak{p}}(a_i \mathfrak{T}_i) = m_{\mathfrak{p}} - 1 \right\}.
\]
Since the minimum defining \(m_{\mathfrak{p}}\) is attained, \(e_{\mathfrak{p}}\ge1\);
on the other hand, \(e_{\mathfrak{p}}\) is at most the number of indices \(i\) with \(a_i\neq0\), because \(v_{\mathfrak{p}}(a_i\mathfrak{T}_i)=m_{\mathfrak{p}}-1\) can hold only when \(a_i\neq0\).

With this notation in place, the general formula of Theorem \ref{mainthm} takes the following explicit product form.
\begin{theorem}\label{the3.3}
Let \(\mathfrak{m} \subseteq R\) be a non-zero ideal, \(\mathfrak{T}_i \mid \mathfrak{m}\) for \(1 \leq i \leq n\), and \(a_i, b\in R\). 
Assume that \(a_i \neq 0\) for at least one \(i\), and consider the restricted linear congruence (\ref{3.1}). 
If there exists a prime ideal \( \mathfrak{p} \mid \mathfrak{m} \) such that either \( m_{\mathfrak{p}} \leq v_{\mathfrak{p}}(\mathfrak{m}) \) and \( v_{\mathfrak{p}}(b) < m_{\mathfrak{p}} - 1 \), or \( m_{\mathfrak{p}} \geq v_{\mathfrak{p}}(\mathfrak{m}) + 1 \) and \( v_{\mathfrak{p}}(b) < v_{\mathfrak{p}}(\mathfrak{m}) \), then the congruence has no solution. 
Otherwise, the number of solutions is given by
\begin{align*}
\prod_{i=1}^n \varphi\left(\frac{\mathfrak{m}}{\mathfrak{T}_i}\right) \prod_{\substack{\mathfrak{p} \mid \mathfrak{m} \\ m_{\mathfrak{p}} \leq v_{\mathfrak{p}}(\mathfrak{m}) \\ v_{\mathfrak{p}}(b) \geq m_{\mathfrak{p}}}} N(\mathfrak{p})^{m_{\mathfrak{p}} - v_{\mathfrak{p}}(\mathfrak{m}) - 1} \left( 1 - \frac{(-1)^{e_{\mathfrak{p}} - 1}}{(N(\mathfrak{p}) - 1)^{e_{\mathfrak{p}} - 1}} \right) \\
\times\prod_{\substack{\mathfrak{p} \mid \mathfrak{m}\\ m_{\mathfrak{p}} \leq v_{\mathfrak{p}}(\mathfrak{m}) \\ v_{\mathfrak{p}}(b) = m_{\mathfrak{p}} - 1}} N(\mathfrak{p})^{m_{\mathfrak{p}} - v_{\mathfrak{p}}(\mathfrak{m}) - 1} \left( 1 - \frac{(-1)^{e_{\mathfrak{p}}}}{(N(\mathfrak{p}) - 1)^{e_{\mathfrak{p}}}} \right),
\end{align*}
where the products are taken over all prime ideals \( \mathfrak{p} \) dividing \( \mathfrak{m}\) satisfying the indicated conditions, and an empty product is interpreted as \(1\). 
\end{theorem}
\begin{proof}
Applying Theorem \ref{mainthm} with \(k=1\), we have \(\mathfrak{D}_i=(a_i)+\mathfrak{m}/\mathfrak{T}_i\), hence \(\mathfrak{T}_i\mathfrak{D}_i=(a_i\mathfrak{T}_i)+\mathfrak{m}\). 
By Theorem
\ref{the2.1},
\[
\eta\left(\frac{\mathfrak{m}}{\mathfrak{D}},\frac{\mathfrak{m}}{\mathfrak{T}_i\mathfrak{D}_i}\right)
=\frac{\varphi\left(\frac{\mathfrak{m}}{\mathfrak{T}_i\mathfrak{D}_i}\right)\mu(w_i)}{\varphi(w_i)},
\]
where, using \(\mathfrak{D}\mid\mathfrak{m}\),
\[
w_i=\frac{\frac{\mathfrak{m}}{\mathfrak{T}_i\mathfrak{D}_i}}{\left(\frac{\mathfrak{m}}{\mathfrak{T}_i\mathfrak{D}_i},\frac{\mathfrak{m}}{\mathfrak{D}}\right)}
=\frac{[\mathfrak{T}_i\mathfrak{D}_i,\mathfrak{D}]}{\mathfrak{T}_i\mathfrak{D}_i}
=\frac{\mathfrak{D}}{(\mathfrak{T}_i\mathfrak{D}_i,\mathfrak{D})}
=\frac{\mathfrak{D}}{((a_i\mathfrak{T}_i,\mathfrak{m}),\mathfrak{D})}
=\frac{\mathfrak{D}}{(a_i\mathfrak{T}_i,\mathfrak{D})}.
\]
Substituting this into Theorem \ref{mainthm}, the number of solutions of (\ref{3.1}) is
\[
\frac{1}{N(\mathfrak{m})} \prod_{i=1}^n \varphi\left( \frac{\mathfrak{m}}{\mathfrak{T}_i} \right)
\sum_{\mathfrak{D}\mid \mathfrak{m}} \eta(b,\mathfrak{D}) \prod_{i=1}^n
\frac{\mu\left( \frac{\mathfrak{D}}{(a_i \mathfrak{T}_i, \mathfrak{D})} \right)}
{\varphi\left( \frac{\mathfrak{D}}{(a_i \mathfrak{T}_i, \mathfrak{D})} \right)}.
\]
By the Chinese remainder theorem, the sum factors over the prime ideals \(\mathfrak{p}\mid\mathfrak{m}\), and it suffices to compute the local sum
\[
W_{\mathfrak{p}} = \sum_{j=0}^{v_{\mathfrak{p}}(\mathfrak{m})}\eta(b,\mathfrak{p}^j)\prod_{i=1}^n
\frac{\mu\left( \frac{\mathfrak{p}^j}{(a_i \mathfrak{T}_i, \mathfrak{p}^j)} \right)}
{\varphi\left( \frac{\mathfrak{p}^j}{(a_i \mathfrak{T}_i, \mathfrak{p}^j)} \right)}
\]
for each \(\mathfrak{p}\mid\mathfrak{m}\), and then multiply the local components.

Assume first that \(m_{\mathfrak{p}}\le v_{\mathfrak{p}}(\mathfrak{m})\). 
If \(0\le j\le m_{\mathfrak{p}}-1\), then \(v_{\mathfrak{p}}(a_i\mathfrak{T}_i)\ge m_{\mathfrak{p}}-1\ge j\) for all \(i\), so \((a_i\mathfrak{T}_i,\mathfrak{p}^j)=\mathfrak{p}^j\) and the product over \(i\) equals \(1\). 
If \(j=m_{\mathfrak{p}}\), then
\((a_i\mathfrak{T}_i,\mathfrak{p}^{m_{\mathfrak{p}}})=\mathfrak{p}^{m_{\mathfrak{p}}-1}\)
for exactly \(e_{\mathfrak{p}}\) indices \(i\), and equals \(\mathfrak{p}^{m_{\mathfrak{p}}}\) for the remaining indices; hence the product equals
\[
\left(\frac{\mu(\mathfrak{p})}{\varphi(\mathfrak{p})}\right)^{e_{\mathfrak{p}}}
=\left(-\frac{1}{N(\mathfrak{p})-1}\right)^{e_{\mathfrak{p}}}.
\]
If \(j>m_{\mathfrak{p}}\), then for some index \(i\) we have
\(\frac{\mathfrak{p}^j}{(a_i\mathfrak{T}_i,\mathfrak{p}^j)}=\mathfrak{p}^{j-(m_{\mathfrak{p}}-1)}\) with exponent at least \(2\), and since \(\mu(\mathfrak{p}^s)=0\) for \(s\ge2\), the entire product vanishes. 
Consequently,
\[
W_{\mathfrak{p}} = \sum_{j=0}^{m_{\mathfrak{p}}-1} \eta(b,\mathfrak{p}^j)
+ \eta(b,\mathfrak{p}^{m_{\mathfrak{p}}})\, \frac{(-1)^{e_{\mathfrak{p}}}}{(N(\mathfrak{p})-1)^{e_{\mathfrak{p}}}}.
\]
By (\ref{2.7}), \(\sum_{j=0}^{r}\eta(b,\mathfrak{p}^j)=N(\mathfrak{p})^r\) if \(r\le v_{\mathfrak{p}}(b)\), and equals \(0\) otherwise; by (\ref{2.8}),
\(\eta(b,\mathfrak{p}^{m_{\mathfrak{p}}})=\varphi(\mathfrak{p}^{m_{\mathfrak{p}}})\)
if \(v_{\mathfrak{p}}(b)\ge m_{\mathfrak{p}}\), and
\(\eta(b,\mathfrak{p}^{m_{\mathfrak{p}}})=-N(\mathfrak{p})^{m_{\mathfrak{p}}-1}\)
if \(v_{\mathfrak{p}}(b)=m_{\mathfrak{p}}-1\). 
Hence
\[
W_{\mathfrak{p}}=
\begin{cases}
N(\mathfrak{p})^{m_{\mathfrak{p}}-1}\left(1 - \dfrac{(-1)^{e_{\mathfrak{p}}-1}}{(N(\mathfrak{p})-1)^{e_{\mathfrak{p}}-1}}\right), & v_{\mathfrak{p}}(b) \ge m_{\mathfrak{p}},\\[6pt]
N(\mathfrak{p})^{m_{\mathfrak{p}}-1}\left(1 - \dfrac{(-1)^{e_{\mathfrak{p}}}}{(N(\mathfrak{p})-1)^{e_{\mathfrak{p}}}}\right), & v_{\mathfrak{p}}(b) = m_{\mathfrak{p}}-1,\\[6pt]
0, & v_{\mathfrak{p}}(b) < m_{\mathfrak{p}}-1.
\end{cases}
\]

If \(m_{\mathfrak{p}} \ge v_{\mathfrak{p}}(\mathfrak{m}) + 1\), then
\(v_{\mathfrak{p}}(a_i \mathfrak{T}_i) \ge m_{\mathfrak{p}}-1\ge v_{\mathfrak{p}}(\mathfrak{m})\) for all \(i\), so the product over \(i\) equals \(1\) for every
\(0\le j\le v_{\mathfrak{p}}(\mathfrak{m})\), and by (\ref{2.7}) we have
\[
W_{\mathfrak{p}} = \sum_{j=0}^{v_{\mathfrak{p}}(\mathfrak{m})} \eta(b,\mathfrak{p}^j) =
\begin{cases}
N(\mathfrak{p})^{v_{\mathfrak{p}}(\mathfrak{m})}, & \text{if } v_{\mathfrak{p}}(b) \ge v_{\mathfrak{p}}(\mathfrak{m}), \\[4pt]
0, & \text{otherwise}.
\end{cases}
\]

Multiplying the local factors and dividing by
\(N(\mathfrak{m})=\prod_{\mathfrak{p}\mid\mathfrak{m}}N(\mathfrak{p})^{v_{\mathfrak{p}}(\mathfrak{m})}\), we obtain the stated product formula. 
This completes the proof.
\end{proof}
We note that the above formula may still evaluate to \(0\) even when neither alternative in the condition above occurs; the complete classification of the cases in which (\ref{3.1}) has no solution is given in Corollary \ref{cor:nosol} below.
\begin{corollary}\label{cor:nosol}
Under the hypotheses of Theorem \ref{the3.3}, the restricted linear congruence (\ref{3.1}) has no solutions if and only if one of the following cases holds:
\begin{enumerate}
\item[(i)] there is a prime ideal \(\mathfrak{p}\mid\mathfrak{m}\) with
\(m_{\mathfrak{p}}\le v_{\mathfrak{p}}(\mathfrak{m})\) and
\(v_{\mathfrak{p}}(b)<m_{\mathfrak{p}}-1\);
\item[(ii)] there is a prime ideal \(\mathfrak{p}\mid\mathfrak{m}\) with
\(m_{\mathfrak{p}}\ge v_{\mathfrak{p}}(\mathfrak{m})+1\) and
\(v_{\mathfrak{p}}(b)<v_{\mathfrak{p}}(\mathfrak{m})\);
\item[(iii)] there is a prime ideal \(\mathfrak{p}\mid\mathfrak{m}\) with
\(m_{\mathfrak{p}}\le v_{\mathfrak{p}}(\mathfrak{m})\),
\(e_{\mathfrak{p}}=1\), and \(v_{\mathfrak{p}}(b)\ge m_{\mathfrak{p}}\);
\item[(iv)] there is a prime ideal \(\mathfrak{p}\mid\mathfrak{m}\) with
\(N(\mathfrak{p})=2\), \(m_{\mathfrak{p}}\le v_{\mathfrak{p}}(\mathfrak{m})\),
\(e_{\mathfrak{p}}\) odd, and \(v_{\mathfrak{p}}(b)\ge m_{\mathfrak{p}}\);
\item[(v)] there is a prime ideal \(\mathfrak{p}\mid\mathfrak{m}\) with
\(N(\mathfrak{p})=2\), \(m_{\mathfrak{p}}\le v_{\mathfrak{p}}(\mathfrak{m})\),
\(e_{\mathfrak{p}}\) even, and \(v_{\mathfrak{p}}(b)=m_{\mathfrak{p}}-1\).
\end{enumerate}
\end{corollary}
\begin{proof}
From the proof of Theorem \ref{the3.3}, the number of solutions is zero exactly when \(W_{\mathfrak{p}}=0\) for some \(\mathfrak{p}\mid\mathfrak{m}\).
By the formulas for \(W_{\mathfrak{p}}\), this is the case in (i) and (ii). 
In the remaining regime \(m_{\mathfrak{p}}\le v_{\mathfrak{p}}(\mathfrak{m})\) and
\(v_{\mathfrak{p}}(b)\ge m_{\mathfrak{p}}-1\), the factor
\[
1-\frac{(-1)^{e_{\mathfrak{p}}-1}}{(N(\mathfrak{p})-1)^{e_{\mathfrak{p}}-1}}
\]
vanishes exactly when \(e_{\mathfrak{p}}=1\), or when
\(N(\mathfrak{p})=2\) and \(e_{\mathfrak{p}}\) is odd; and the factor
\[
1-\frac{(-1)^{e_{\mathfrak{p}}}}{(N(\mathfrak{p})-1)^{e_{\mathfrak{p}}}}
\]
vanishes exactly when \(N(\mathfrak{p})=2\) and \(e_{\mathfrak{p}}\) is
even, since
\(\left|\frac{(-1)^{e_{\mathfrak{p}}}}{(N(\mathfrak{p})-1)^{e_{\mathfrak{p}}}}\right|<1\)
for \(N(\mathfrak{p})\ge3\). This gives (iii), (iv) and (v).
\end{proof}
The following special case is worth recording separately. 
When \(m_{\mathfrak{p}}=1\) for every \(\mathfrak{p}\mid\mathfrak{m}\), the product formula simplifies considerably and generalizes the Rademacher--Brauer formula (\ref{1.3}).
\begin{corollary}\label{cor3.3}
Assume that for every prime ideal \(\mathfrak{p}\mid\mathfrak{m}\) one has \(m_{\mathfrak{p}}=1\), that is, for each \(\mathfrak{p}\mid\mathfrak{m}\) there exists at least one index \(i\in\{1,\dots,n\}\) with \(v_{\mathfrak{p}}(a_i\mathfrak{T}_i)=0\); 
in particular, \(a_i\neq0\) for at least one \(i\), so Theorem \ref{the3.3} applies. 
Then the number of solutions of the restricted linear congruence (\ref{3.1}) is
\[
\frac{1}{N(\mathfrak{m})} \prod_{i=1}^n \varphi\left( \frac{\mathfrak{m}}{\mathfrak{T}_i} \right)
\prod_{\substack{\mathfrak{p} \mid \mathfrak{m} \\ v_{\mathfrak{p}}(b) \geq 1}}
\left( 1 - \frac{(-1)^{e_{\mathfrak{p}} - 1}}{(N(\mathfrak{p}) - 1)^{e_{\mathfrak{p}} - 1}} \right)
\prod_{\substack{\mathfrak{p} \mid \mathfrak{m} \\ v_{\mathfrak{p}}(b) = 0}}
\left( 1 - \frac{(-1)^{e_{\mathfrak{p}}}}{(N(\mathfrak{p}) - 1)^{e_{\mathfrak{p}}}} \right).
\]
\end{corollary}
We close this section with two explicit examples. 
Both are chosen so that the solution sets can be determined by hand, and the formulas of Theorem \ref{mainthm} and the preceding corollaries are checked against these direct counts. 
The first takes place in the ring of integers of a number field whose class number is not one, the second in a complete local ring outside the global setting.
\begin{example}\label{ex:Qsqrt-5}
Let \(R=\mathbb{Z}[\sqrt{-5}]\), the ring of integers of \(\mathbb{Q}(\sqrt{-5})\), whose class number is \(2\). 
Here \((2)=\mathfrak{p}_2^2\) with the non-principal prime ideal \(\mathfrak{p}_2=(2,1+\sqrt{-5})\) of norm \(2\), and the reduced residue classes modulo \((2)\) are \(1\) and \(\sqrt{-5}\).

Consider the restricted congruence
\[
x_1+x_2\equiv 1+\sqrt{-5}\pmod{(2)},\qquad (x_j)+(2)=(1)\quad (j=1,2).
\]
The only solutions are \((1,\sqrt{-5})\) and \((\sqrt{-5},1)\), so the number of solutions is \(2\).

In Theorem \ref{mainthm} take \(k=1\), \(n=2\), \(a_1=a_2=1\), \(b=1+\sqrt{-5}\), \(\mathfrak{T}_1=\mathfrak{T}_2=(1)\), so that \(\mathfrak{D}_1=\mathfrak{D}_2=(1)\) and the Euler-function factor equals \(1\). 
Since \((b)=\mathfrak{p}_2\mathfrak{p}_3\), where \(\mathfrak{p}_3=(3,1+\sqrt{-5})\), we have \((b,(2))=\mathfrak{p}_2\), and Theorem \ref{the2.1} gives
\[
\eta(b,(1))=1,\quad \eta(b,\mathfrak{p}_2)=1,\quad \eta(b,(2))=-2,
\]
together with
\[
\eta((2),(2))=2,\quad \eta(\mathfrak{p}_2,(2))=-2,\quad
\eta((1),(2))=\mu((2))=0,
\]
the last one because \((2)=\mathfrak{p}_2^2\) is not square-free. 
Hence
\[
\frac{1}{N((2))}\sum_{\mathfrak{D}\mid(2)}\eta(b,\mathfrak{D})\,
\eta\!\left(\frac{(2)}{\mathfrak{D}},(2)\right)^{2}
=\frac{1}{4}(4+4+0)=2,
\]
in agreement with the direct count. 
Alternatively, Corollary \ref{cor3.3} applies with \(m_{\mathfrak{p}_2}=1\) and \(e_{\mathfrak{p}_2}=2\), giving
\[
\frac{1}{N((2))}\,\varphi((2))^2\left(1+\frac{1}{N(\mathfrak{p}_2)-1}\right)
=\frac{4}{4}\cdot 2=2.
\]
\end{example}
\begin{example}\label{ex:Z3}
Let \(R=\mathbb{Z}_3\), the ring of \(3\)-adic integers. 
It is a complete discrete valuation ring, hence a Dedekind domain with the finite norm property; its non-zero ideals are the powers \((3^r)\), with \(N((3^r))=3^r\). 
The reduced residue classes modulo \((9)\) are \(1,2,4,5,7,8\).

Consider the restricted congruence
\[
x_1+x_2+x_3\equiv 0\pmod{(9)},\qquad
(x_j)+(9)=(1)\ (j=1,2),\quad (x_3)+(9)=(3),
\]
so that \(x_1,x_2\) are units modulo \((9)\), while \(x_3\in\{3,6\}\). 
If \(x_3=3\), then \(x_1+x_2\equiv6\pmod9\), and for each unit \(u\) modulo \(9\) the element \(6-u\equiv-u\pmod3\) is again a unit, so there are \(6\) pairs. 
The same argument applies when \(x_3=6\), with \(x_1+x_2\equiv3\pmod9\). 
Hence the number of solutions is \(12\).

In Theorem \ref{mainthm} take \(k=1\), \(n=3\), \(a_1=a_2=a_3=1\), \(b=0\), \(\mathfrak{T}_1=\mathfrak{T}_2=(1)\), \(\mathfrak{T}_3=(3)\), so that \(\mathfrak{D}_1=\mathfrak{D}_2=\mathfrak{D}_3=(1)\) and the Euler-function factor equals \(1\). 
The divisors of \((9)\) are \((1)\), \((3)\), \((9)\), and Theorem \ref{the2.1} gives
\[
\eta(0,(1))=1,\quad \eta(0,(3))=2,\quad \eta(0,(9))=6,
\]
together with
\[
\eta((9),(9))=6,\quad \eta((9),(3))=2,\qquad
\eta((3),(9))=-3,\quad \eta((3),(3))=2,
\]
while \(\eta((1),(9))=\mu((9))=0\). 
Hence
\[
\frac{1}{N((9))}\sum_{\mathfrak{D}\mid(9)}\eta(0,\mathfrak{D})\,
\eta\!\left(\frac{(9)}{\mathfrak{D}},(9)\right)^{2}
\eta\!\left(\frac{(9)}{\mathfrak{D}},(3)\right)
=\frac{1}{9}(72+36)=12,
\]
in agreement with the direct count. 
Alternatively, Theorem \ref{the3.3} applies with \(m_{(3)}=1\), \(e_{(3)}=2\), and \(v_3(0)=\infty\ge m_{(3)}\), giving
\[
\varphi((9))^{2}\,\varphi((3))\cdot N((3))^{1-2-1}
\left(1-\frac{(-1)^{2-1}}{(N((3))-1)^{2-1}}\right)
=72\cdot\frac{1}{9}\cdot\frac{3}{2}=12.
\]
\end{example}
\begin{acknowledgment}
The author would like to thank Associate Professor Su Hu for proposing the topic of this paper, which has been a source of inspiration throughout this work, and Dr.~Tianfang Qi for carefully reading the entire manuscript and offering valuable suggestions on the writing and the structure.
\end{acknowledgment}


\begin{thebibliography}{99}
\bibitem{Alder1958}
H.~L. Alder, A generalization of the Euler $\phi$-function, Amer. Math. Monthly {\bf 65} (1958), 690--692.

\bibitem{Babu2024}
C.~G. Karthick~Babu, R. Bera and B. Sury, Linear congruences in several variables with congruence restrictions, Monatsh. Math. {\bf 211} (2026), no.~1, 1--24.

\bibitem{Bibak2017}
K.~Bibak, B.~Kapron, V.~Srinivasan, R.~Tauraso, and L.~T\'oth,
Restricted linear congruences, J.~Number Theory {\bf 171} (2017), 128--144.

\bibitem{Brauer1926}
A.~Brauer, L\"osung der Aufgabe 30, Jahresber. Dtsch. Math.-Ver. {\bf 35} (1926), 92--94.

\bibitem{Cohen1949}
E. Cohen, An extension of Ramanujan's sum, Duke Math. J. {\bf 16} (1949), 85--90.

\bibitem{Cohen1955}
E. Cohen, A class of arithmetical functions, Proc. Nat. Acad. Sci. U.S.A. {\bf 41} (1955), 939--944.

\bibitem{Cohen1959}
E. Cohen, Representations of even functions $({\rm mod}$ $r)$. III. Special topics, Duke Math. J. {\bf 26} (1959), 491--500.

\bibitem{Deaconescu2000}
M. Deaconescu, Adding units mod $n$, Elem. Math. {\bf 55} (2000), no.~3, 123--127.

\bibitem{Dixon1960}
J.~D. Dixon, A finite analogue of the Goldbach problem, Canad. Math. Bull. {\bf 3} (1960), 121--126.

\bibitem{Jacobson1972}
D.~N. Jacobson and K.~S. Williams, On the number of distinguished representations of a group element, Duke Math. J. {\bf 39} (1972), 521--527.

\bibitem{Lehmer1913}
D.~N. Lehmer, Certain theorems in the theory of quadratic residues, Amer. Math. Monthly {\bf 20} (1913), 151--157.

\bibitem{Namboothiri2018}
K.~V. Namboothiri, On the number of solutions of a restricted linear congruence, J. Number Theory {\bf 188} (2018), 324--334.

\bibitem{Narkiewicz2004}
W. Narkiewicz, {\it Elementary and analytic theory of algebraic numbers}, third edition, 
Springer Monographs in Mathematics, Springer, Berlin, 2004.

\bibitem{Nicol1954}
C. Nicol and H.~S. Vandiver, A von Sterneck arithmetical function and restricted partitions with respect to a modulus, Proc. Nat. Acad. Sci. U.S.A. {\bf 40} (1954), 825--835.

\bibitem{Qi2025}
T. Qi, A generalized Ramanujan sum over a residually finite Dedekind domain, Acta Math. Hungar. {\bf 175} (2025), no.~2, 333--351.

\bibitem{Rademacher1925}
H.~Rademacher, Aufgabe 30, Jahresber. Dtsch. Math.-Ver. {\bf 34} (1925), 158.

\bibitem{Ramanujan1918}
S.~A. Ramanujan, On certain trigonometrical sums and their applications in the theory of numbers, Trans. Cambridge Philos. Soc. {\bf 22} (1918), 259--276.

\bibitem{Rearick1963}
D. Rearick, A linear congruence with side conditions, Amer. Math. Monthly {\bf 70} (1963), 837--840.

\bibitem{Sander2009}
J. Sander, On the addition of units and nonunits mod $m$, J. Number Theory {\bf 129} (2009), no.~10, 2260--2266.

\bibitem{Sander2013}
J. Sander and T. Sander, Adding generators in cyclic groups, J. Number Theory {\bf 133} (2013), no.~2, 705--718.

\bibitem{Sburlati2003}
G. Sburlati, Counting the number of solutions of linear congruences, Rocky Mountain J. Math. {\bf 33} (2003), no.~4, 1487--1497.

\bibitem{Sun2014}
C.~F. Sun and Q.-H. Yang, On the sumset of atoms in cyclic groups, Int. J. Number Theory {\bf 10} (2014), no.~6, 1355--1363.

\bibitem{Zariski1965}
O. Zariski and P. Samuel, {\it Commutative algebra. Vol. 1}, corrected reprinting of the 1958 edition, 
Graduate Texts in Mathematics, No. 28, Springer, New York-Heidelberg-Berlin, 1975.

\bibitem{zheng2022}
Z. Zheng, M. Chen and J. Xu, On Gauss sums over Dedekind domains, Int. J. Number Theory {\bf 18} (2022), no.~5, 1047--1063.

\bibitem{Zheng2023}
Z. Zheng, M. Chen and Z. Hong, On Ramanujan sums over a Dedekind domain with finite norm property, Acta Math. Sin. (Engl. Ser.) {\bf 39} (2023), no.~1, 149--160.
\end{thebibliography}
\end{document}